\documentclass[11pt]{article}

\usepackage{amsmath,amssymb,latexsym,theorem,enumerate,QED}
\usepackage[centering,dvips]{geometry}
\usepackage{hyperref}
\usepackage[all]{xypic}
\def\xybiglabels{\def\labelstyle{\textstyle}}

{\theorembodyfont{\rmfamily}\newtheorem{example}{Example}[section]}
{\theorembodyfont{\rmfamily}\newtheorem{remark}[example]{Remark}}
{\theorembodyfont{\rmfamily}\newtheorem{Def}[example]{Definition}}
{\theorembodyfont{\rmfamily}}
{\theorembodyfont{\rmfamily}}
\newtheorem{prop}[example]{Proposition}
\newtheorem{cor}[example]{Corollary}
\newtheorem{thm}[example]{Theorem}

\newcommand{\labto}[1]{\stackrel{#1}{\longrightarrow}}

\def\xybiglabels{\def\labelstyle{\textstyle}}

\def\Set{\mathsf{Set}}

\def\qed{\hfill $\Box$}

\def\rho{\varrho}
\def\clap#1{\hbox to 0pt{hss#1\hss}}

\def\Gpds{\mathsf{Gpd}}

\def\Top{\mathsf{Top}}

\def\Z{\mathbb{Z}}

\def\leq{\leqslant}
\def\geq{\geqslant}
\def\Groups{\mathsf{Groups}}

\newcommand{\directs}[2]{\def\objectstyle{\scriptstyle} \objectmargin={0pt}
\xy
(0,4)*+{}="a",(0,-2)*+{\rule{0em}{1.5ex}#2}="b",(7,4)*+{\;#1}="c"
\ar@{->} "a";"b" \ar @{->}"a";"c" \endxy }
\newcommand{\xdirects}[2]{\def\objectstyle{\scriptstyle}
\objectmargin={0pt} \xy
(0,0)*+{}="a",(0,-6)*+{\rule{0em}{1.5ex}#2}="b",(7,0)*+{\;#1}="c"
\ar@{->} "a";"b" \ar @{->}"a";"c" \endxy}

\newcommand{\threeaxes}[3]{\def\objectstyle{\scriptstyle}  \objectmargin={0pt}
\xy
(0,0)*+{}="a",(0,-6)*+{\rule{0em}{1.5ex}#2}="b",(7,0)*+{\;#1}="c",
(14,-3)*+{\;#3}="d" \ar@{->} "a";"b" \ar @{->}"a";"c"  \ar
@{->}"a";"d"\endxy }

\newcommand{\threeaxesb}[3]{\def\objectstyle{\scriptstyle}  \objectmargin={0pt}
\xy
(0,0)*+{}="a",(0,-8)*+{\rule{0em}{1.5ex}#2}="b",(10,0)*+{\;#1}="c",
(6,9)*+{\;#3}="d" \ar@{->} "a";"b" \ar @{->}"a";"c"  \ar
@{->}"a";"d"\endxy }

\def\Ob{\operatorname{Ob}}

\def\rel{\rm{rel}}
\def\tab{^{\rm totab}}
\def\ab{^{\rm ab}}

\begin{document}
\title{Covering morphisms of groupoids, derived modules and a
1-dimensional Relative Hurewicz Theorem\thanks{MSCClassification:
20L99,20J99,55P99; KEYWORDS: Hurewicz theorems, fundamental group,
derived module, covering spaces and covering morphisms of
groupoids}}
\author{Ronald Brown\thanks{School of Computer Science, Bangor University, Gwynedd LL57 1UT, UK; r.brown@bangor.ac.uk}}
\maketitle
\begin{abstract}
We fill  a lacuna in the literature by giving a version in dimension
1 of the Relative Hurewicz Theorem, and relate this to
abelianisations of groupoids, covering spaces,   covering morphisms
of groupoids, and Crowell's notion of derived modules.
\end{abstract}

\section*{Introduction}
A classical result is that when the space $X$ with base point $a$ is
path connected then the  canonical morphism
$$\pi_1(X,a)\to H_1(X)$$ is just abelianisation of the fundamental group $\pi_1(X,a)$ of $X$ at  $a$.

However there is a notion of the fundamental groupoid $\pi_1(X,A)$
on {\it a set $A$ of base points}, introduced in \cite{B67} to allow
for  an associated Seifert--van Kampen Theorem for non connected
spaces, and this concept is clearly also necessitated in discussing
spaces with a group of operators, as developed in
\cite{Higgins-Taylor} and \cite[Chapter 11]{B2006}.

Our main aim in this paper is to give circumstances which enable us
to identify the canonical morphism
$$\pi_1(X,A) \to H_1(X,A)$$
as an abelianisation of the  groupoid $\pi_1(X,A)$.  This
identification is clearly a 1-dimensional form of the Relative
Hurewicz Theorem. We then use covering morphisms of groupoids to
relate this result to Crowell's notion of derived module,
\cite{Cr71}, and show its relevance to covering spaces.

The usual Relative Hurewicz Theorem (RHT) states that if $n \geq 2$,
and $(X,A,x)$ is an $(n-1)$-connected pair of connected spaces, then
the natural Hurewicz morphism in dimension $n$
$$\pi_n(X,A,x) \to H_n(X,A)$$
is given by factoring the action of $\pi_1(A,x)$. This result is
basic in algebraic topology. It has been rephrased in \cite[Section
8]{BH81:col}, see also \cite{BS:fibredcats},  as the statement that
under the same conditions the natural morphism
$$\pi_n(X,A,x) \to \pi_n(X \cup CA,CA,x)$$has the same properties:
 the proofs of these results are given in \cite{BH81:col,BHS}  independently of
homology theory, as a special case of a Higher Homotopy Seifert-van
Kampen Theorem, involving pushouts, or more general colimits. This
more general theorem even has a version involving a set $X_0$ of
base points, for which one has to use the notion of induced module
or crossed module over groupoids for a precise statement of a
theorem analogous to the RHT, see for example \cite{BS:fibredcats}.

\section{Derived modules}\label{sec:derivedmodules}
In this section we recall the work of Crowell in \cite{Cr71}. Let
$\phi: F \to G$ be a morphism of groups, and let $M$ be a (right)
$G$-module.
\begin{Def}A {\it $\phi$-derivation} $d: F \to M$ is a function
satisfying
$$d(uv)= d(u) ^ {\phi v} + d(v)$$
for all $u,v \in F$.

A $\phi$-derivation $\partial: F \to D$ is said to be a {\it
universal $\phi$-derivation} if it is a $\phi$-derivation and for
any other $\phi$-derivation $d: F \to M$ there is a unique
$G$-morphism $d': D \to M$ such that $d' \partial=d$.

In case $\phi$ is the identity morphism then these terms are
abbreviated by omitting the $\phi$. An example of a derivation is
$e: G \to IG$, $g \mapsto g-1$,  where $IG$ is the augmentation
ideal of $G$. \qed
\end{Def}
\begin{example}
An example of a universal derivation is $e: G \to IG$ where $IG$ is
the augmentation ideal, which is generated as $G$-module by the
elements $g-1, g \ne 1, g \in G$ and $e$ is given by $g \mapsto
g-1$. \qed
\end{example}
Since a universal $\phi$-derivation $\partial:F \to D $ determines
$D$ uniquely up to isomorphism, the $G$-module $D$ is then called
the {\it derived module} of the morphism $\phi$, and is written
$D_\phi$. The derived module may be constructed directly and also as
$IF \otimes_{\Z F}\Z G$.

One of the key uses of the derived module is in the following exact
sequence of Crowell, \cite[Section 4]{Cr71}.
\begin{prop}\label{prop:exactseq}
Let $\phi: F \to G$ be an epimorphism of groups, with kernel $N$.
Then there is an exact sequence of $G$-modules
\begin{equation}
0 \to N \ab \to D_\phi \to IG \to 0.
\end{equation}
\end{prop}
\begin{remark}
This construction of an exact sequence of modules from an exact
sequence of groups is the key to relating the classification of
extensions of groups to the homological algebra of modules, and is
treated in \cite{Mac63}. The construction of the derived module is
generalised to the groupoid case in \cite{brohig5}, as is essential
for the results there, and the above  sequence is also generalised.
\qed
\end{remark}

\section{Abelianisations of groupoids}
\label{ssec:II7groupoidabelian}

We use at several points the well known abelianisation of a group.
This construction gives a functor
$$^{ab} \colon \Groups \to \mathsf{Ab}$$
which is left adjoint to the inclusion of categories
$$\mathsf{Ab} \to \Groups.$$

We define a groupoid $G$ to be {\it abelian} if all its vertex
groups are abelian groups. The {\it abelianisation of a groupoid}
$G$, which we write $G \ab$, is obtained by quotienting $G$ with the
normal subgroupoid generated by the commutators from all vertex
groups.

Nonetheless, we need another kind of abelianisation that associates
to each groupoid not an abelian groupoid but an abelian group and a
morphism $\upsilon \colon  G \to G^{\tab}$ which is universal for
morphisms to abelian groups. We call $G \tab$ the {\it universal
abelianisation } of the groupoid $G$.

Let $\mathsf{Ab}$, $\Groups$, $\Gpds$ denote respectively the
categories of abelian groups, groups, and groupoids. Each of the
inclusions
\begin{equation}
\mathsf{Ab} \to \Groups \to \Gpds
\end{equation}
has a left adjoint. That from groupoids to groups is called the {\it
universal group} $UG$ of a groupoid $G$ and is described in detail
in \cite{H71} and \cite[Section 8.1]{B2006}. In particular, the
universal group of a groupoid $G$ is the free product of the
universal groups of the transitive components of $G$. Any transitive
groupoid $G$ may be written in a non canonical way as the free
product $G(a_0) * T$ of a vertex group $G(a_0)$ and an indiscrete or
tree groupoid $T$. Then
$$UG \cong G(a_0)* UT$$ and $UT$ is the free group on the elements
$x \colon  a_0 \to a$ in $T$ for all $a \in \Ob(T)$, $a \ne a_0$.

It follows that the universal abelianisation is given by
\begin{equation}G \tab \cong (UG)\ab,
\end{equation}and
also that  $G \tab$ is isomorphic to the direct sum of the $G_i\tab$
over all components $G_i$ of $G$. So for a transitive groupoid $G$
with $ a_0 \in \Ob G$
\begin{equation}\label{eq:totabgpd}
G \tab \cong G(a_0) \ab \oplus F\end{equation} where $F$ is the free
abelian group on the elements $x \colon  a_0 \to a$ in $T$ for all
$a \in \Ob(T), a\ne a_0$, for $T$ a wide tree subgroupoid of $G$.

If $G$ is a totally disconnected groupoid on the set $X$, and
$\theta  \colon  G \to X$ is the unique morphism over $X$ to the
discrete groupoid on $X$, then the derived module of $\theta $ is
the above defined abelianisation $G^{ab}$ of $G$.

\section{Covering morphisms of groupoids} \label{ssecII5:covGpdCrs}
For the convenience of readers, and to fix the notation, we recall
here the basic facts on covering morphisms of groupoids.  The
earliest definition of covering morphism of groupoids seems to be in
\cite{PASmith1,PASmith2}, where such a morphism is called a
`regular' morphism.   The ideas were developed independently in
\cite{H64}. Many basic facts which should be seen as a part of
`combinatorial groupoid theory'  can be found in the books
\cite{H71,B2006}. However we find it convenient to adopt different
conventions from, say, \cite{B2006}, focussing on costars rather
than stars.

Let $G$ be a groupoid. For each object $a_0$ of $G$ the {\it Costar}
of $a_0$ in $G$, denoted by ${\rm Cost}_{G}~a_0$, is the union of
the sets $G(a,a_0)$ for all objects $a$ of $G$, i.e. ${\rm
Cost}_{G}a_0 = \{g \in G \mid tg = a_0\}$. A morphism $p  \colon
\widetilde{G} \rightarrow G$ of groupoids is a {\it covering
morphism} if for each object $\widetilde{a}$ of $\widetilde{G}$ the
restriction of $p$ \begin{equation} \label{eq:appcostars}   {\rm
Cost}_{\widetilde{G}}~\widetilde{a} \to {\rm
Cost}_{G}~p\widetilde{a} \end{equation} is bijective. In this case
$\widetilde{G}$ is called a {\it covering groupoid of G}.
\label{GSApA:Gd14} More generally $p$ is called a {\it fibration of
groupoids}, \cite{B70},  if the restrictions of $p$ to the Costars
as in Equation \eqref{eq:appcostars} are surjective.
\label{app:deffibgpds}

A basic result for covering groupoids is { \em unique path lifting}.
That is, let $p  \colon  \widetilde{G} \rightarrow G$ be a covering
morphism of groupoids, and let $(g_1, g_2, \ldots, g_n)$ be a
sequence of composable elements of $G$. Let $\tilde{a} \in
Ob(\widetilde{G})$ be such that $p\tilde{a} $ is the target  of
$g_n$. Then there is a unique composable sequence $(\tilde{g}_1,
\tilde{g}_2, \ldots, \tilde{g}_n) $ of elements of $\widetilde{G}$
such that $\tilde{g}_n$ ends at $\tilde{a}$ and $p\tilde{g}_i = g_i,
i=1, \ldots, n$.

If $G$ is a groupoid, the slice category ${\Gpds\mathsf{Cov}}/G$
\label{GSApA:Gd15} of coverings of $G$ has as objects the covering
morphisms $p  \colon  H \rightarrow G$ and has as morphisms the
commutative diagrams of morphisms of groupoids, where $p$ and $q$
are covering morphisms,
$$\xybiglabels\xymatrix{H \ar[rr]^f \ar [rd]_p & & K \ar [ld]^q  \\
& K & }
$$

By a standard result on  compositions and covering morphisms
(\cite[10.2.3]{B2006}), $f$ also is a covering morphism. It is
convenient to write such a diagram as a triple $(f,p,q)$. The
composition in ${\Gpds\mathsf{Cov}}/G$ is then given as usual by
$$(g,q,r)(f,p,q) = (gf,p,r). $$

It is a standard result (see for example \cite{H71,B70}) that the
category ${\Gpds\mathsf{Cov}}/G$ is equivalent to the category of
operations of the groupoid $G$ on sets. We give the definitions and
notations which we will use for this equivalence.

Recall we are writing composition of $g \colon  p\to q$ and $h
\colon q \to r$ in a groupoid as $gh \colon  p \to r$. This is the
opposite of the notation for functions in the category $\Set$; the
composite of a function  $f \colon X \to Y$ and $g \colon  Y \to Z$
is $gf \colon  X \to Z$ with value $(gf)(x)=g(f(x))$. It is possible
to resolve this confusion by writing functions on the right of their
argument as $(x)f$: This `algebraist's' convention is followed
successfully in \cite{H71}, and contrasts with the usual `analyst's'
convention. Because of this `opposite' nature of our conventions we
have to make the following definition.

\begin{Def}
A {\it left operation of a groupoid $G$ on sets} is a functor $X
\colon  G^{\rm op} \to \Set$. If $p \in \Ob(G)$, $g \colon  p \to q$
in $G$, and $x \in X(p)$, then $X(g)(x)\in X(q)$ may also be written
$g \cdot x$.
\end{Def}

Thus if $X \colon  G^{\rm op} \to \Set$ is a functor, then the {\it
action groupoid} $\widetilde{G}=G \ltimes X$ has object set the
disjoint union of the sets $X(p)$ for $p \in Ob(G)$ and morphisms $x
\to y$ the pairs $(g,x)$ such that $x \in X(tg)$ and $y= X(g)x$; in
operator notation: $(g,x) \colon  gx \to x$. The composition is
$(g',gx)(g,x)=(g'g,x)$. The projection morphism $G \ltimes X \to G,
\; (g,x) \mapsto g$, is  a covering morphism.

This `semidirect product' or `Grothendieck construction' is
fundamental for constructing covering morphisms to the groupoid $G$.
This so \index{Grothendieck construction}called `Grothendieck
construction' has also been developed by C. Ehresmann in
\cite{Ehr-Gatt}, in which he defines both an action of a category
and the associated `category of hypermorphisms',
\index{hypermorphism}and also what in the case of local groupoids he
calls the complete enlargement of a species of structures. For
example, if $a_0$ is an object of the transitive groupoid $G$, and
$A$ is a subgroup of the object group $G(a_0)$ then the groupoid $G$
operates on the family of cosets $\{gA \mid g \in {\rm Cost}_G~a_0
\}$, by $g'.(gA) =g'gA$ whenever $g'g$ is defined, and the
associated covering morphism $\widetilde{G} \to G $ defines the
covering groupoid $\widetilde{G}$ of the groupoid $G$ determined by
the subgroup $A$. When $A$ is trivial this gives the {\it universal
cover} at $a_0$ of the groupoid $G$.

In particular, this gives the universal covering groupoid of a group
$G$, whose objects are the elements of $G$ and morphisms are pairs
$(g,h) \colon  gh \to h$, $g,h \in G$, and composition is
$$(g,hk)(h,k)=(gh,k)$$for $g,h,k \in G$.  Then $G$ operates on the right of the
universal cover by $(g,h) ^k=(g,hk)$, for $g,h,k \in G$. This
operation preserves the map $p$ and is called a {\it covering
transformation}. The relevance of these ideas to derivations is as
follows. Let $d:G \to \widetilde{G}$ be defined by $g\mapsto (g,1)$.
Then $d$ satisfies $d(gh)=d(g)^h \, d(h)$ as illustrated by the
following diagram:

$$\def\labelstyle{\textstyle}
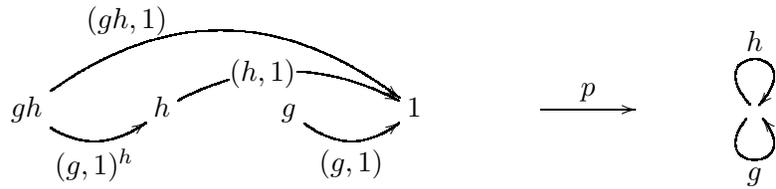
\xymatrix@C=3pc{gh\ar @/_1pc/[r]_{(g,1)^h}\ar@/^2.5pc/[rrr]^(0.3){(gh,1)} &
h \ar@/^1.2pc/[rr]|(0.4){(h,1)}
  &\ar@/_1pc/ [r]_{(g,1)}g& 1& \ar[r]^p &&\ar@(ul,ur)[]^h \ar@(dl,dr)[]_g  }$$

\begin{example} 
Here is a simple example: the universal covering groupoid
$\widetilde{K}$ of the Klein 4-group $K= \mathsf{C}_2 \times
\mathsf{C}_2$ with elements say $1,a,b,ab$. This group is generated
by $a,b$  with the relations $a^2,b^2,aba^{-1}b^{-1}$, which we
write respectively $r,s,t$. Then $\widetilde{K}$ has the elements of
$K$ as vertices and a morphism $(g,x) \colon gx \to x$ for each $g,x
\in K$. The covering morphism $p \colon \widetilde{K} \to K$ is
$(g,x) \mapsto g$. In terms of the generators $a,b$ we obtain a
diagram of $\widetilde{K}$  as the left hand diagram in the
following picture:
$$\xybiglabels 
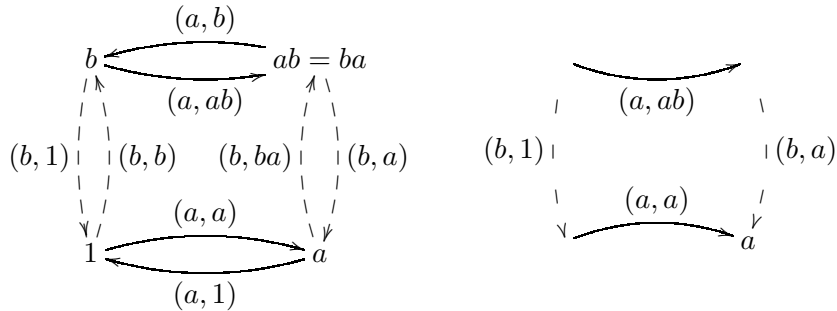
\xymatrix@=5pc{b \ar@/_1ex/ @{-->}[d]_{(b,1)} \ar @/_1.5ex/ [r]_{(a,ab)}& ab=ba
\ar@/^1.5ex/ @{-->}[d]^{(b,a)}\ar @/_1.5ex/ [l]_{(a,b)}\\
1 \ar@/_1.5ex/ @{-->}[u]_{(b,b)}\ar@/^1.5ex/ [r]^{(a,a)}& a
\ar@/^1.5ex/ @{-->}[u]^{(b,ba)}\ar@/^1.5ex/ [l]^{(a,1)}}\qquad
\xymatrix@=5pc{\ar@/_1ex/ @{-->}[d]_{(b,1)} \ar @/_1.5ex/
[r]_{(a,ab)}&\ar@/^1.5ex/ @{-->}[d]^{(b,a)}\\\ar@/^1.5ex/
[r]^{(a,a)} &a} $$ Note that for example $(a,ab) \colon b \to ab$
because $a^2=1$. The right hand diagram illustrates a lift of the
loop $b^{-1}a^{-1}ba$ in $K$ to a loop starting and ending at $a$ in
the diagram of $\widetilde{K}$. You should draw the similar loops
starting in turn at $1,b,ab$. It is also possible to view  these
four loops as  forming boundaries of four `lifts' of the relation
$t$. \qed
\end{example}

For a space $X$ which is locally path-connected and semi-locally
$1$-connected there is a standard analysis of covering spaces of $X$
in terms of operations of the fundamental groupoid $\pi_1X$ on sets,
and this is the formalism usually followed in texts on algebraic
topology; however there is another formulation of this result, which
goes back in the simplicial set case to \cite{GZ67}, and states that
under the above local conditions on $X$ the fundamental groupoid
functor
\begin{equation}\label{eq:equivpi1}
\pi_1:{\Top\mathsf{Cov}}/X \to {\Gpds\mathsf{Cov}}/\pi_1
X\end{equation} is an equivalence of categories. The construction of
the inverse equivalence in Equation \eqref{eq:equivpi1} is of course
closely related to the usual method of constructing covering spaces,
but seems to fall more naturally as part of this theorem, since it
relates morphisms in two different contexts, namely spaces and
groupoids, and the proofs are base point free.  We indicate how the
inverse equivalence may be constructed.

Let then $q: \widetilde{G} \to \pi_1 X $ be a covering morphism of
groupoids. Let $\widetilde{X}= \Ob(\widetilde{G}) $ and let
$p=\Ob(q): \widetilde{X} \to X$. One uses the properties of covering
morphisms of groupoids and the assumed local conditions on $X$ to
lift sufficiently small neighbourhoods of points of  $X$ to
neighbourhoods of points of $\widetilde{X}$ to obtain a base for a
topology on $\widetilde{X}$ with the required properties. Full
details of this construction are in \cite[Sections 10.5,10.6
]{B2006}, and indeed were given in the 1968 edition of that book.

The above equivalence is used in \cite{BMu} to discuss covering
groups of non-connected topological groups.

\section{Covering morphisms and derived modules}
In this section we link the notions of covering morphisms with those
of derived modules, in effect clarifying and giving algebraic
versions of results in \cite[Section 5]{Cr71} and \cite[Section
11]{W49:CHII}.

Let $\phi: F \to G$ be a morphism of groups with kernel $N$. We form
the universal covering groupoid $p \colon  \widetilde{G} \to G$ and
the pullback
\begin{equation}
\xybiglabels \vcenter{\xymatrix{\widehat{F} \ar [d] _q \ar [r]
^{\bar{\phi}}& \widetilde{G} \ar [d] ^p\\
F \ar [r] _\phi & G}}
\end{equation}
Note that a morphism in $\widehat{F}$ is a pair $(u,(\phi u,g))
\colon (\phi u)g \to g$, $u \in F, g \in G$. Since $u$ determines
$\phi u$, we can write a morphism  of $\widehat{F}$ as $(u,g) \colon
(\phi u)g \to g$. One utility of this diagram is that if $F$ is the
free group on a set $S$ of generators, and $\phi$ is an epimorphism,
then the graph $q^{-1}(S)$ can be seen as the Cayley graph of $G$
for this set of generators, and $\widehat{F}$ is the free groupoid
on this graph. The free groupoid on a graph was introduced in
\cite{H64} and a recent use is in \cite{crisp-paris-artinconj}.

The group $G$ operates on the right of $\widetilde{G}$ by
$(h,g)^k=(h,gk)$, $k \in G$ and so on the right of $\widehat{F}$ by
$(u,g)^k=(u,gk)$, $k \in G$.

It is easy to prove, and alternatively follows from the
Mayer-Vietoris sequence of the pullback, \cite[10.7.6]{B2006} and
\cite{BHK-MayV}, that $\widehat{F}$ is connected if and only if
$\phi$ is surjective, and that $q$ maps the vertex group
$\widehat{F}(1)$ bijectively to $N$.

\begin{thm}\label{thm:exactsequences}
If $\phi: F \to G$ is an epimorphism of groups with kernel $N$ then
there  is a natural isomorphism of exact sequences of $G$-modules
\begin{equation}
\vcenter{\xymatrix{0 \ar[r] &  N \ab _{\phantom{x}} \ar [r] \ar [d]&
D_\phi ^{\phantom{X}} \ar [r]\ar [d]^\eta & IG
\ar [r]\ar [d] & 0\\
0 \ar [r] & N \ab _{\phantom{x}}\ar [r] & \widehat{F}\tab
_{\phantom{x}}\ar [r] & \widetilde{G}\tab  _{\phantom{x}}\ar [r] &0
.}}
\end{equation}
\end{thm}
\begin{proof}
The top line of the diagram is the standard exact sequence of
modules derived from the exact sequence of groups $1 \to N
\labto{\phi} F \to G \to 1$ as shown in \cite[Section 5]{Cr71} and
which  corresponds to a classical derivation of such a module
sequence, \cite[Section IV.6]{Mac63}.

Here $D_\phi$ is defined by a universal $\phi$-derivation $d: F \to
D_\phi$. Thus to construct $\eta$ we need a $\phi$-derivation $e: F
\to \widehat{F} \tab$. This is the composition $F \labto{f}
\widehat{F} \labto{a} \widehat{F} \tab$, where $f$ is the function
$u \mapsto (u,1)$. Then $$f(uv)=(uv,1)= (uv,v)(v,1) =(u,1)^{\phi
v}(v,1)$$ from which it follows that $e$ is a $\phi$-derivation.

Now the vertex group $\widehat{F}(1)$ is isomorphic to $N$ under $p$
so we identify these. Also $\Ob(\widetilde{G})= G $ and so, since
$\widehat{F}$ is connected, $$\widehat{F}\tab \cong N \ab \oplus (G
\times G)\tab.$$

However  $(G \times G)\tab$ is the free abelian group on elements $g
\in G, g \ne 1$, by Equation \eqref{eq:totabgpd}. However in keeping
with the module action of $G$ this can also be regarded as the free
abelian group on elements $g-1$, $g \ne 1$ where
$$hg-1=(h-1)^g +(g-1).$$
This is another description of the augmentation module $IG$.

Thus we have the map of exact sequences as given in the Theorem, and
so by the 5-lemma, $\eta$ is an isomorphism.
\end{proof}

\section{Homology and homotopy } \label{sec:IIIhomologyandhomotopy}
The homology groups of a cubical set $K$  are defined as follows.
First we form the chain complex $C'(K)$ where $C'_n(K)$ is
\index{homology!cubical}the free abelian group on $K_n$, and with
boundary
\begin{equation}
  \partial k= \sum _{i=1}^n (-1)^i (\partial^-_i k - \partial ^+_i
  k).
\end{equation}
It is easily verified that this gives a chain complex, i.e.
$\partial \partial =0$. However if $K$ is a point, i.e. $K_n$ is a
singleton for all $n$, then the homology groups of $C'(K)$  are
$\mathbb Z$ in even dimensions, whereas we want the homology of a
point to be zero in dimensions $> 0$. We therefore \index{normalise
cubical chains} normalise, i.e. factor $C'(K)$ by the subchain
complex  generated by the degenerate cubes. This gives the chain
complex $C_*(K)$ of $K$, and the homology groups of this chain
complex are defined to be the homology groups of $K$. In particular
the (cubical) singular homology groups of the space $X$ are defined
to be the homology groups of the cubical singular complex $S^\square
X$.

A full exposition of this cubical homology theory is in
\cite{hilton-wylie,Mas80}. It is proved in
\cite{Eil-MacLane-acyclic} using acyclic models that the cubical
singular homology groups are isomorphic to the simplicial singular
homology groups.  Recent works using cubical methods are
\cite{desc-cub,blanc-johnson-turner,isaacson-cubsets}.

Let $X_*$ be a filtered space: that is, $X_*$ is given by a space
$X$ and a sequence of subspaces
$$ X_*  \colon = \quad  X_0 \subseteq X_1 \subseteq X_2 \subseteq \cdots \subseteq
X_n \subseteq \cdots \subseteq X $$ which we call a {\it filtration}
of $X$. A map of filtered spaces $f: X_* \to Y_*$ is a map $f: X \to
Y$ preserving the filtration.

\begin{example}
\begin{enumerate}[1.]
\item A standard example of a filtered space is a $CW$-complex with its
skeletal filtration. So among these we have the $n$-cubes $I^n_*$
with their usual cell structure and skeletal filtration.
\item Let $A$ be a subspace of the topological space $X$, and let $n
\geq 0$. Then we have a filtration of $X$ which we could write
$X^n_A$ defined by $(X^n_A)_i= A$ for $i \leq n$ and $(X^n_A)_i=X$
for $i
>n$. These filtrations occur in proofs of the Hurewicz Theorem.
\item Another example is the free monoid $JX$ on a
topological space $X$, with $JX$ filtered by the length of the
words. And finally, a filtration on a manifold is determined by a
Morse function on the manifold.  \qed
\end{enumerate}
\end{example}

The notion of map of filtered spaces  allows us to define the {\it
filtered singular complex} $RX_*$ which in dimension $n$ is the set
of filtered maps $I^n_* \to X_*$. This has the structure of cubical
set, and other structure such as connections and compositions, which
are discussed in \cite{BH81:algcub}. Also $RX_*$ is a fibrant
cubical set, or in other terms,  is a Kan complex, as is easily
proved: the use of this structure as a foundation for algebraic
topology is described in \cite{BHS}.

The cubical chain complex $C_*(X_*)$ is defined to be the normalised
chain complex of the cubical set $RX_*$.

The following retraction Proposition, which is taken from
\cite[Section 9]{BH81:col}, is one step towards the Hurewicz
theorem. This retraction proposition should be compared with a
special case discussed in \cite[Section III.7]{Mas80}. The history
of classical papers on singular homology and the Hurewicz Theorem
shows the use of deformation theorems of the type of this
Proposition, as for example in Blakers \cite{Bl48}. However the use
of cubical methods as here, rather than of simplicial methods and
chain complexes, does seem to simplify the proof, partly since
cubical methods are easier for constructing homotopies.

In the proof, a useful lemma is that if $(Y,Z)$ is a cofibred pair,
and $f \colon (Y,Z) \to (X,A)$ is a map of pairs which is deformable
(as a map of pairs) into $A$, then $f$ is deformable into $A$ rel
$Z$ (\cite[7.4.4]{B2006}).

\begin{prop} \label{III2prop:RtoKhomequ}Let $X_*$ be a filtered
space such that the following conditions $\psi (X_*, m)$ hold for
all $m \geq 0$:
\begin{enumerate}[]
       \item $\psi (X_*, 0) :$ The map $\pi_0 X_0 \rightarrow
  \pi_0 X$ induced by inclusion is surjective;
      \item $\psi (X_*, 1) :$ Any path in $X$ joining points of $X_0$ is deformable in $X$ rel
 end points to a path in $X_1$; \item $\psi (X_*, m) (m \geq 2
 ) :$ For all $\nu \in X_0$ , the map
$$\pi_m (X_m , X_{m-1} , \nu ) \rightarrow \pi_m (X, X_{m-1} , \nu
)$$ induced by inclusion is surjective. \end{enumerate} Then the
inclusion $i  \colon  RX_* \rightarrow KX=S^\square X$ is a homotopy
equivalence of cubical sets. \end{prop}
\begin{proof}   There exist maps
$h_m  \colon  K_m X \rightarrow K_{m+1} X,\; r_m  \colon  K_m X
\rightarrow K_m X$ for $m \geq 0$ such that
\begin{enumerate}[(i)]
     \item $\partial^-_{m+1} h_m = 1, \partial^+_{m+1} h_m = r_m $;
    \item   $r_m (KX) \subseteq R_m X_*$
  and $h_m \mid R_m X_* = \varepsilon_{m+1} $;
    \item $\partial^\tau_i h_m = h_{m-1} \partial^\tau_i$ for $1 \leq i
\leq m$ and $\tau = -, +$;
    \item $h_m \varepsilon_j = \varepsilon_j h_{m-1}$ for
    $1 \leq j \leq m$.
\end{enumerate}
Such $r_m , h_m$ are easily constructed by induction, starting with
$h_{-1} = \emptyset$ , and using $\psi (X_*, m)$ to define $h_m
\alpha$ for elements $\alpha$ of $K_m X$ which are not degenerate
and do not lie in $R_m X_*.$ Here is a picture for $h_1$:
\begin{equation*}
\xybiglabels \xymatrix@M=0pt@=5pc{\bullet
\ar@{}[dr]|{h_1k}\ar@{--}[d]_(-0.07){r\partial^-_1 k}
_(1.03){r\partial^+_1 k}_{rk} \ar@{-} [r] ^{h_0\partial^-_1k} &
\bullet \ar@{-} [d] ^(-0.07){\partial ^-_1k} ^(1.03){\partial ^+_1k}^k \\
\bullet \ar@{-} [r] _{h_0
\partial^+_1k} & \bullet  }
\end{equation*}

These maps define a retraction $r \colon KX \rightarrow RX_*$ and a
homotopy $h \simeq ir$ rel $RX_*$. \end{proof}
\begin{cor}\label{III2cor:RtoKisoforCW} If $X_*$ is the skeletal
filtration of a CW-complex, then the inclusion $RX_* \to S^\square
X$ is a homotopy equivalence of fibrant  cubical sets. \qed
\end{cor}

\begin{cor}\label{III2cor:RtoKiso}If the conditions
$\psi (X_*, m)$ of the proposition hold for all $m \geq 0,$ then the
inclusion $i \colon RX_* \rightarrow S^\square X$ induces a homotopy
equivalence of chain complexes and hence an isomorphism of all
homology and homotopy groups.
\end{cor}
\begin{proof}
The result on homotopy is standard, and that on homology follows
from the development in \cite{Mas80}.
\end{proof}

\subsection{Relative Hurewicz Theorem: dimension 1}
In this subsection we identify the total abelianisation of the
groupoid $\pi_1(X,A)$,  see Section \ref{ssec:II7groupoidabelian},
in certain cases in terms of homology.

\begin{Def} We now coin a term: for a subspace
$A$ of $X$, let $C_*(X \rel_0 A )$ denote the sub chain complex of
$C_*(X^0_A)$ in which the only change is that $C_0(X \rel_0 A)=0$;
thus all elements of $C_1(X \rel_0 A)$ are cycles, and of course all
generating elements of $C_n(X \rel_0 A)$ map vertices of $I^n$ into
$A$.

We write $H_*(X \rel_0 A)$ for the homology of this chain complex.
\qed
\end{Def}
\begin{thm} \label{II1:abelianisationfundgpd}Let $A$ be a subspace of the space $X$. Then a Hurewicz
morphism $$\omega \colon  \pi_1(X,A) \to H_1(X\rel_0 A)$$ is defined
and induces an  isomorphism $$\omega' \colon \pi_1(X,A) \tab \to
H_1(X \rel_0 A). $$
\end{thm}
\begin{proof}
For each path class $[f] \in \pi_1(X,A)$ the representative $f$
determines a generator of $C_1(X\rel_0 A)$. Differing choices of $f$
yield homologous elements of $C_1(X\rel_0 A)$, so this defines
$\omega$ as a function. If $f \circ g$ is a composite of paths with
vertices in $A$  then the diagram
\begin{equation}
\xybiglabels \vcenter{\xymatrix@M=0pt@=3pc{\ar @{-} [d] _f|\tip  \ar
@{-} [r]|\tip ^{f \circ g}& \ar @{-} [d]|\tip ^1 \\ \ar @{-}
[r]_g|\tip &}}
\end{equation}
extends to a map of $I^2 \to X$ with vertices mapped to $A$ whose
boundary shows that $\omega$ is a morphism to $H_1(X\rel_0 A)$. It
hence defines $\omega' \colon  \pi_1(X,A) \tab \to H_1( X\rel_0 A)$.

Now $C_1(X\rel_0 A)$ is free abelian on the non degenerate paths $f
\colon I \to X$ with vertices in $A$. So a morphism  $\eta \colon
C_1( X\rel_0 A) \to \pi_1(X,A) \ab$ is defined by sending $f$ to its
class in $\pi_1(X,A) \ab$. It is easy to check that $\eta
\partial_2=0$, so that $\eta$ defines a morphism $$H_1(X \rel_0 A)
\to \pi_1(X,A)\tab, $$ and that $\eta$ is inverse to $\omega'$.
\end{proof}

Next we relate $H_*(X \rel_0 A)$ to the standard relative homology.
The result we need generalises the classical case,  when $X$ is path
connected and $A$ consists of a single point,\cite[Lemma
7.2]{Mas80}.
\begin{prop}\label{prop:III2chneq}
If $A$ meets each path component of $X$, then the inclusion
\[C_*(X^0_A) \to C_*(X)\] is a chain equivalence.
\end{prop}
\begin{proof} This is an immediate consequence of Corollary
\ref{III2cor:RtoKiso}.
\end{proof}

We say $C_*(A)$ is {\it concentrated in dimension $0$} if $C_i(A)=0$
for $i>0$. This occurs for example if $A$ is  totally path
disconnected, and so if $A$ is discrete. \index{Hurewicz!Relative
Theorem!in dimension 1}
\begin{thm}[Relative Hurewicz Theorem: dimension
1]\label{thm:III2relHurdim1} If $A$ is totally path disconnected and
meets each path component of $X$  then the natural morphism
\[  \pi_1(X,A)^{\rm totab} \to  H_1(X,A) \]
is an isomorphism.
\end{thm}
\begin{proof}
We define $A_*$ to be the constant filtered space with value $A$. So
we regard $A_*$ as a sub-filtered space of $X^0_A$.

We consider the morphism of exact sequences of chain complexes

\begin{equation}
\xybiglabels \vcenter{\xymatrix{0 \ar [r] & C_*(A_*)\ar [d]^= \ar [r] & C_*(X^0_A) \ar [r] \ar [d]^i& C_*(X^0_A,A_*)
\ar [r] \ar [d]^j& 0\\
0 \ar [r] & C_*(A) \ar [r]
& C_*(X) \ar [r]   &
C_*(X,A) \ar [r]& 0 \\
 }}
\end{equation}
where classically the lower sequence defines relative homology
$H_*(X,A)$, and the upper sequence defines $H_*(X^0_A,A_*)$. Under
our assumptions, the morphism $i$ is a homotopy equivalence and
hence so also is $j$ (since all the chain complexes are free in each
dimension).

Our assumption that $A$ is totally path disconnected implies that
$C_i(A)=0$ for $i>0$. This implies that $C_*(X^0_A,A_*) \cong C_*(X
\rel_0 A)$. So the Theorem follows from  Theorem
\ref{II1:abelianisationfundgpd} and Proposition
\ref{prop:III2chneq}.
\end{proof}

\begin{thm}
Let $p: \widehat{X} \to X$ be a covering map of connected spaces
determined by the normal subgroup $N$ of $F= \pi_1(X,x)$, and let
$\phi: F \to G=F/N$ be the quotient morphism. Let
$\widetilde{X}_0=p^{-1}(x)$ so that there is a bijection
$$\widetilde{X}_0 \cong G = F/N.$$
Let $\widehat{F}$ be as in Theorem \ref{thm:exactsequences}. Then
$\widehat{F} \cong \pi_1(\widehat{X},\widetilde{X}_0)$ and so there
is an isomorphism
$$\widehat{F}\tab \cong H_1(\widehat{X},\widetilde{X}_0).$$
Hence there is a universal $\phi$-derivation $\pi_1(X,x) \to
H_1(\widehat{X},\widetilde{X}_0)$.
\end{thm}
\begin{proof}The covering morphism of groupoids
$\pi_1(\widehat{X},\widetilde{X}_0) \to F= \pi_1(X,x)$ is clearly
isomorphic to $\widehat{F} \to F$ since they are both determined by
the normal subgroup $N$ of $F$. But now Theorem
\ref{thm:III2relHurdim1} gives the result.
\end{proof}
\begin{remark} This description of $\pi_1(X,x) \to
H_1(\widehat{X},\widetilde{X}_0)$ as a universal $\phi$-derivation
is essentially the result of Crowell \cite[Section 5]{Cr71}. It also
gives the case $i=1$ of \cite[Proposition 5.2]{brohig5} which
relates the fundamental crossed complex of certain filtered spaces
to the chains with operators defined by universal covers. That
Proposition generalises  work in \cite[Section 11]{W49:CHII};
Whitehead's use of  the chains of the universal cover in this and a
number of papers was strongly influenced by Reidemeister's paper
\cite{Rei1934}. That work developed  also into work of
Eilenberg--Mac~Lane on homology of spaces with operators,
\cite{eilenberg-operatorsI,eil-macl-operatorsII}, as well as into
Fox's work on the free differential calculus \cite{Fox53}. \qed
\end{remark}

\end{document}